\documentclass[10pt,fleqn]{article}

\usepackage[includeheadfoot]{geometry}
\usepackage{times}
\usepackage{graphics}
\usepackage{graphicx}
\usepackage{amsbsy}
\usepackage{amssymb}
\usepackage{amsfonts}
\usepackage{amstext}
\usepackage{amsmath}
\usepackage{booktabs}
\usepackage{multirow}
\usepackage{dsfont}
\usepackage{psfrag}
\usepackage[parfill]{parskip}

\usepackage{fancyhdr}
\fancypagestyle{plain}{%
\fancyhead[L]{}
\fancyhead[C]{}
\fancyhead[R]{}
\fancyfoot[L]{\small\sf Davide Giglio}
\fancyfoot[C]{\small\sf DIBRIS -- University of Genova}
\fancyfoot[R]{\sf\thepage}

}

\usepackage{amsthm}
\newtheoremstyle{thstyle}{6pt}{3pt}{}{}{\bf}{.}{.5em}{}

\theoremstyle{thstyle}

\title{{\sc\large TECHNICAL REPORT}\\\vspace{25pt}A MILP model for single machine family scheduling\\with sequence-dependent batch setup\\and controllable processing times}
\author{Davide Giglio\vspace{6pt}\\
{\small Department of Informatics, Bioengineering, Robotics and Systems Engineering (DIBRIS)}\vspace{-2pt}\\
{\small University of Genova}\vspace{-2pt}\\
{\small Via Opera Pia 13, 16145 -- Genova, Italy}\vspace{-2pt}\\
{\small davide.giglio@unige.it}%
\vspace{20pt}}
\date{\includegraphics[scale=.2]{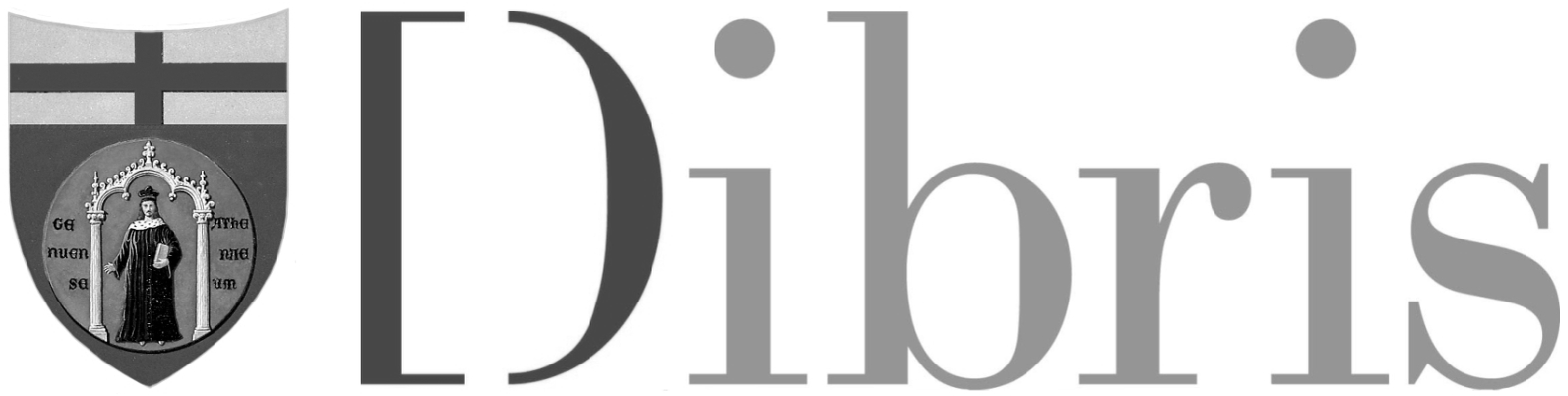}}

\begin{document}

\setlength{\belowdisplayskip}{9pt}
\setlength{\belowdisplayshortskip}{0pt}
\setlength{\abovedisplayskip}{9pt}
\setlength{\abovedisplayshortskip}{0pt}

\pagestyle{fancy}

\maketitle

\vspace{20pt}

\begin{abstract}
A mathematical programming model for a class of single machine family scheduling problems is described in this technical report, with the aim of comparing the performance in solving the scheduling problem by means of mathematical programming with the performance obtained when using optimal control strategies, that can be derived from the application of a dynamic programming-based methodology proposed by the Author. The scheduling problem is characterized by the presence of sequence-dependent batch setup and controllable processing times; moreover, the generalized due-date model is adopted in the problem. Three mixed-integer linear programming (MILP) models are proposed. The best one, from the performance point of view, is a model which makes use of two sets of binary variables: the former to define the relative position of jobs and the latter to define the exact sequence of jobs. In addition, one of the model exploits a stage-based state space representation which can be adopted to define the dynamics of the system.
\end{abstract}

\newpage

\section{Introduction}

A mathematical programming model for the single machine family scheduling problem described in the paper {\itshape ``Optimal control strategies for single machine family scheduling with sequence-dependent batch setup and controllable processing times''} \cite{GiglioJOSH} is here presented in order to compare the dynamic programming-based approach proposed in \cite{GiglioJOSH} with a different methodology. Among the possible formulations, a mixed-integer linear programming (MILP) model which makes use of two sets of binary variables (one to define the relative position of jobs and the other to define the exact sequence of jobs) is adopted; as a matter of fact, a preliminary performance analysis has shown that such a model outperforms other models which solve the same problem.

This technical report is organized as follows. The scheduling problem is described in section~\ref{sec:scheduling} whereas the MILP model is reported in section~\ref{sec:model1}. The test of the model (on the numerical example considered in \cite{GiglioJOSH}) and an experimental analysis (on randomly-generated instances of different sizes of the problem) are in section~\ref{sec:test} and~\ref{sec:experiment}, respectively. A comparison of the times required to find an optimal solution between the algorithm proposed in \cite{GiglioJOSH} (implemented and solved with Matlab) and the MILP model (implemented and solved with Cplex) is in section~\ref{sec:comparison}. Two alternative MILP formulations are reported in section~\ref{sec:model23} (the second of them exploits the stage-based state space representation adopted in \cite{GiglioJOSH}).

\section{The scheduling problem} \label{sec:scheduling}

Consider a single machine where $N_{k}$ jobs of class $P_{k}$, $k = 1, \ldots, K$, $K > 1$, have to be executed. All jobs belonging to the same class are equivalent. All jobs are available at time instant $0$, and preemption is not allowed. A sequence of $N_{k}$ due-dates for jobs of class $P_{k}$, namely $dd_{k,1}, dd_{k,2}, \ldots,$ $dd_{k,N_{k}}$, $k = 1, \ldots, K$, is specified. It is assumed that $dd_{k,i} \leq dd_{k,i+1}$ for any $k = 1, \ldots, K$ and $i = 1, \ldots, N_{k}-1$, and that jobs are assigned to the due-dates according to the EDD (earliest due-date) rule, in accordance with the generalized due-date model.

The processing (or service) time of the $i$-th job of class $P_{k}$ is a continuous variable $pt_{k,i} \in [pt^{\mathrm{low}}_{k} , pt^{\mathrm{nom}}_{k}]$, $k = 1, \ldots, K$, $i = 1, \ldots, N_{k}$; it is assumed that, once the processing time of a job has been chosen, it cannot be changed during the service. $pt^{\mathrm{nom}}_{k}$ is the nominal value of the processing time, which can be compressed up to $pt^{\mathrm{low}}_{k}$. In accordance with a resource consumption model, the processing time can be expressed as $pt_{k,i} (u_{k,i}) = pt^{\mathrm{nom}}_{k} - \gamma_{k} \, u_{k,i}$, with $0 \leq u_{k,i} \leq \overline{u}_{k}$, $\overline{u}_{k} = (pt^{\mathrm{nom}}_{k} - pt^{\mathrm{low}}_{k})/\gamma_{k}$, where $u_{k,i}$ is the amount of a continuous resource used to reduce the processing time of the $i$-th job of class $P_{k}$, and $\gamma_{k}$ is the positive compression rate for jobs of class $P_{k}$.

A costly setup is required between the execution of two jobs of different classes. In this connection, let $st_{k,h}$ (resp., $sc_{k,h}$), $k = 1, \ldots, K$, $h = 1, \ldots, K$, be the setup time (cost) which has to be spent (paid) when switching from a job of class $P_{k}$ to one of class $P_{h}$, $h \ne k$. No setup is required between the execution of two jobs of the same class (i.e., $st_{k,k}=0$ and $sc_{k,k}=0$ $\forall \, k = 1,\ldots,K$). Besides, the execution of the first job does not require any setup, regardless of the class. Then, $st_{0,k} = 0$ and $sc_{0,k} = 0$, $k = 1,\ldots,K$, having denoted by $st_{0,k}$ and $sc_{0,k}$ the initial setup time and the initial setup cost, respectively, when a job of class $P_{k}$ is the first processed job.

The performance criterion (cost function) to be minimized is the sum of the total weighted tardiness, the total weighted consumption cost (of the continuous resource used to compress the processing time), and the total setup cost, namely
\begin{equation} \label{equ:performance_criterion}
\sum_{k=1}^{K} \sum_{i=1}^{N_{k}} \alpha_{k,i} \, T_{k,i} + \sum_{k=1}^{K} \sum_{i=1}^{N_{k}} \beta_{k} \, \gamma_{k} \, u_{k,i} + \sum_{k=1}^{K} \sum_{i=1}^{N_{k}} sc_{\pi_{k,i},k}
\end{equation}
where: $\alpha_{k,i}$ is the unitary tardiness cost for the $i$-th job of class $P_{k}$; $T_{k,i}$ is the tardiness of the $i$-th job of class $P_{k}$; $\beta_{k}$ is the unitary cost related to the deviation from the nominal service time for the execution of jobs of class $P_{k}$; $\pi_{k,i}$ is the index of the class of the job which precedes the $i$-th job of class $P_{k}$ in the job sequence, or is the value $0$ if the $i$-th job of class $P_{k}$ is the first processed job. The set of values $\Pi = \{Ê\pi_{k,i}$ , $k = 1, \ldots, K$, $i = 1, \ldots, N_{k} \}$ is sufficient to define a specific job sequence; moreover, let $U = \{Êu_{k,i}$ , $k = 1, \ldots, K$, $i = 1, \ldots, N_{k} \}$. The total weighted tardiness depends on both $\Pi$ and $U$, as $T_{k,i} = \max \{ S_{k,i}(\Pi,U) + st_{\pi_{k,i},k} + pt^{\mathrm{nom}}_{k} - \gamma_{k} \, u_{k,i} - dd_{k,i} , 0 \}$, being $S_{k,i}(\Pi,U)$ the time instant at which the $i$-th job of class $P_{k}$ starts its execution when values $\pi_{k,i}$ and $u_{k,i}$, $k = 1, \ldots, K$, $i = 1, \ldots, N_{k}$ are adopted.

The objective is to find a job sequence $\Pi$ and the values $u_{k,i}$, $k = 1, \ldots, K$, $i = 1, \ldots, N_{k}$, which minimize the cost function~\eqref{equ:performance_criterion}.

\section{The mathematical programming model} \label{sec:model1}

\subsection{Variables and parameters} \label{sec:varpar}

\subsubsection{Main decision variables}

\begin{itemize}
\item $x_{h,j,k,i}$ : binary variable which defines the relative position of jobs
\begin{itemize}
\item $x_{h,j,k,i} = 1$ if the $j$-th job of class $P_{h}$ is executed before the $i$-th job of class $P_{k}$
\item $x_{h,j,k,i} = 0$ otherwise
\end{itemize}
\item $\delta_{h,j,k,i}$ : binary variable which defines the sequence of jobs
\begin{itemize}
\item $\delta_{h,j,k,i} = 1$ if the $i$-th job of class $P_{k}$ is the next job after the execution of the $j$-th job of class $P_{h}$
\item $\delta_{h,j,k,i} = 0$ otherwise
\end{itemize}
\item $u_{k,i}$ : amount of a continuous used to reduce the processing time of the $i$-th job of class $P_{k}$
\end{itemize}

\subsubsection{Other decision variables}

\begin{itemize}
\item $S_{k,i}$ : time instant at which the $i$-th job of class $P_{k}$ starts its execution
\item $pt_{k,i}$ : processing time of the $i$-th job of class $P_{k}$
\item $T_{k,i}$ : tardiness of the $i$-th job of class $P_{k}$
\item $\Omega_{k,i}$ : setup cost which is paid for the $i$-th job of class $P_{k}$
\item $\Lambda_{k,i}$ : setup time which is spent for the $i$-th job of class $P_{k}$
\end{itemize}

\subsubsection{Parameters}

\begin{itemize}
\item $N_{k}$ : number of jobs of class $P_{k}$
\item $pt^{\mathrm{nom}}_{k}$ : nominal value of the processing time of jobs of class $P_{k}$
\item $pt^{\mathrm{low}}_{k}$ : lowest value of the processing time of jobs of class $P_{k}$
\item $dd_{k,i}$ : due-date for the $i$-th job of class $P_{k}$
\item $st_{h,k}$ : setup time which is spent when switching from a job of class $P_{h}$ to a job of class $P_{k}$
\item $sc_{h,k}$ : setup cost which is paid when switching from a job of class $P_{h}$ to a job of class $P_{k}$
\item $\alpha_{k,i}$ : unitary tardiness cost for the $i$-th job of class $P_{k}$
\item $\beta_{k}$ : unitary processing time's deviation cost for jobs of class $P_{k}$
\item $\gamma_{k}$ : processing time's compression rate for jobs of class $P_{k}$
\item $M$ : large and positive value (``big number'' larger than an upper-bound of the maximum completion time)
\end{itemize}

\subsection{The MILP model} \label{sec:MIPmodel}

\subsubsection{Objective function}

\begin{equation}
\min \quad \sum_{k=1}^{K} \sum_{i=1}^{N_{k}} \alpha_{k,i} \, T_{k,i} + \sum_{k=1}^{K} \sum_{i=1}^{N_{k}} \beta_{k} \, \gamma_{k} \, u_{k,i} + \sum_{k=1}^{K} \sum_{i=1}^{N_{k}} \Omega_{k,i}
\end{equation}

\subsubsection{Constraints}

\begin{equation} \label{cst:tardiness}
T_{k,i} \geq S_{k,i} + \Lambda_{k,i} + pt_{k,i} - dd_{k,i} \qquad \forall k = 1, \ldots, K \quad \forall i = 1, \ldots, N_{k}
\end{equation}
\begin{equation} \label{cst:setupcost}
\Omega_{k,i} = \sum_{h=1}^{K} \sum_{j=1}^{N_{h}} sc_{h,k} \; \delta_{h,j,k,i} \qquad \forall k = 1, \ldots, K \quad \forall i = 1, \ldots, N_{k}
\end{equation}
\begin{equation} \label{cst:setuptime}
\Lambda_{k,i} = \sum_{h=1}^{K} \sum_{j=1}^{N_{h}} st_{h,k} \; \delta_{h,j,k,i} \qquad \forall k = 1, \ldots, K \quad \forall i = 1, \ldots, N_{k}
\end{equation}
\begin{equation} \label{cst:alljobs}
\sum_{h=1}^{K} \sum_{j=1}^{N_{h}} \sum_{k=1}^{K} \sum_{i=1}^{N_{k}} \delta_{h,j,k,i} = \bigg( \sum_{k=1}^{K} N_{k} \bigg) - 1
\end{equation}
\begin{equation} \label{cst:onlyonepredecessor}
\sum_{h=1}^{K} \sum_{j=1}^{N_{h}} \delta_{h,j,k,i} \leq 1 \qquad \forall k = 1, \ldots, K \quad \forall i = 1, \ldots, N_{k}
\end{equation}
\begin{equation} \label{cst:onlyonesuccessor}
\sum_{k=1}^{K} \sum_{i=1}^{N_{k}} \delta_{h,j,k,i} \leq 1 \qquad \forall h = 1, \ldots, K \quad \forall j = 1, \ldots, N_{h}
\end{equation}
\begin{equation} \label{cst:genduedatedelta1}
\sum_{j=1}^{i} \delta_{k,i,k,j} = 0 \qquad \forall k = 1, \ldots, K \quad \forall i = 1, \ldots, N_{k}
\end{equation}
\begin{equation} \label{cst:genduedatedelta2}
\sum_{j=i+2}^{N_{k}} \delta_{k,i,k,j} = 0 \qquad \forall k = 1, \ldots, K \quad \forall i = 1, \ldots, N_{k}-2
\end{equation}
\begin{equation} \label{cst:consumption}
u_{k,i} \leq \dfrac{pt^{\mathrm{nom}}_{k} - pt^{\mathrm{low}}_{k}}{\gamma_{k}} \qquad \forall k = 1, \ldots, K \quad \forall i = 1, \ldots, N_{k}
\end{equation}
\begin{equation} \label{cst:processingtime}
pt_{k,i} = pt^{\mathrm{nom}}_{k} - \gamma_{k} \, u_{k,i} \qquad \forall k = 1, \ldots, K \quad \forall i = 1, \ldots, N_{k}
\end{equation}
\begin{multline} \label{cst:startafter}
S_{k,i} \geq S_{h,j} + \Lambda_{h,j} + pt_{h,j} - M \, ( 1 - x_{h,j,k,i}) \\ \forall h = 1, \ldots, K \quad \forall j = 1, \ldots, N_{h} \quad \forall k = 1, \ldots, K \quad \forall i = 1, \ldots, N_{k} \quad (h,j) \neq (k,i)
\end{multline}
\begin{multline} \label{cst:startbefore}
S_{h,j} \geq S_{k,i} + \Lambda_{k,i} + pt_{k,i} - M \, x_{h,j,k,i} \\ \forall h = 1, \ldots, K \quad \forall j = 1, \ldots, N_{h} \quad \forall k = 1, \ldots, K \quad \forall i = 1, \ldots, N_{k} \quad (h,j) \neq (k,i)
\end{multline}
\begin{equation} \label{cst:genduedatex1}
x_{k,j,k,i} = 1 \qquad \forall k = 1, \ldots, K \quad \forall i = 1, \ldots, N_{k} \quad \forall j = 1, \ldots, i-1
\end{equation}
\begin{equation} \label{cst:genduedatex2}
x_{k,j,k,i} = 0 \qquad \forall k = 1, \ldots, K \quad \forall i = 1, \ldots, N_{k} \quad \forall j = i, \ldots, N_{k}
\end{equation}
\begin{multline} \label{cst:xnocirculartwo}
x_{h,j,k,i} + x_{k,i,h,j} = 1 \\ \forall h = 1, \ldots, K \quad \forall j = 1, \ldots, N_{h} \quad \forall k = 1, \ldots, K \quad \forall i = 1, \ldots, N_{k} \quad (h,j) \neq (k,i)
\end{multline}
\begin{multline} \label{cst:xnocircularthree}
x_{h,j,k,i} + x_{k,i,l,m} + x_{l,m,h,j} \leq 2 \\
\begin{array}{r} \forall h = 1, \ldots, K \quad \forall j = 1, \ldots, N_{h} \quad \forall k = 1, \ldots, K \quad \forall i = 1, \ldots, N_{k} \\ (h,j) \neq (k,i) \quad (k,i) \neq (l,m) \quad (l,m) \neq (h,j) \end{array}
\end{multline}
\begin{multline} \label{cst:relationshipdeltax}
x_{h,j,k,i} \geq 1 - M \, ( 1 - \delta_{h,j,k,i}) \\ \forall h = 1, \ldots, K \quad \forall j = 1, \ldots, N_{h} \quad \forall k = 1, \ldots, K \quad \forall i = 1, \ldots, N_{k}
\end{multline}
\begin{equation} \label{cst:defx}
x_{h,j,k,i} \in \{ 0 , 1 \} \qquad \forall h = 1, \ldots, K \quad \forall j = 1, \ldots, N_{h} \quad \forall k = 1, \ldots, K \quad \forall i = 1, \ldots, N_{k}
\end{equation}
\begin{equation} \label{cst:defdelta}
\delta_{h,j,k,i} \in \{ 0 , 1 \} \qquad \forall h = 1, \ldots, K \quad \forall j = 1, \ldots, N_{h} \quad \forall k = 1, \ldots, K \quad \forall i = 1, \ldots, N_{k}
\end{equation}
\begin{equation} \label{cst:defu}
u_{k,i} \geq 0 \qquad \forall k = 1, \ldots, K \quad \forall i = 1, \ldots, N_{k}
\end{equation}
\begin{equation} \label{cst:defS}
S_{k,i} \geq 0 \qquad \forall k = 1, \ldots, K \quad \forall i = 1, \ldots, N_{k}
\end{equation}
\begin{equation} \label{cst:defpt}
pt_{k,i} \geq 0 \qquad \forall k = 1, \ldots, K \quad \forall i = 1, \ldots, N_{k}
\end{equation}
\begin{equation} \label{cst:defT}
T_{k,i} \geq 0 \qquad \forall k = 1, \ldots, K \quad \forall i = 1, \ldots, N_{k}
\end{equation}
\begin{equation} \label{cst:defOmega}
\Omega_{k,i} \geq 0 \qquad \forall k = 1, \ldots, K \quad \forall i = 1, \ldots, N_{k}
\end{equation}
\begin{equation} \label{cst:defLambda}
\Lambda_{k,i} \geq 0 \qquad \forall k = 1, \ldots, K \quad \forall i = 1, \ldots, N_{k}
\end{equation}

\vspace{12pt}
Constraint~\eqref{cst:tardiness} computes, taking into account~\eqref{cst:defT} and the presence of $T_{k,i}$ in the cost function, the tardiness of the $i$-th job of class $P_{k}$ which is executed. Constraints~\eqref{cst:setupcost} and~\eqref{cst:setuptime} compute, respectively, the setup cost and the setup time which are paid/spent when executing the $i$-th job of class $P_{k}$. Constraint~\eqref{cst:alljobs} ensures that all jobs are executed. Constraint~\eqref{cst:onlyonepredecessor} (resp., \eqref{cst:onlyonesuccessor}) guarantees that each job can have one predecessor (successor) only; the inequality is necessary to deal with the first (the last) executed job. Constraints~\eqref{cst:genduedatedelta1} and~\eqref{cst:genduedatedelta2} apply the generalized due-date model: with respect of class $P_{k}$, only jobs $i$ and  $i+1$ can be (but not necessarily are) consecutive jobs in the overall sequence of job executions. Constraint~\eqref{cst:consumption} determines, taking into account~\eqref{cst:defu}, the interval which is allowed for the amount of a continuous resource used to reduce the processing time of the $i$-th job of class $P_{k}$. Constraint~\eqref{cst:processingtime} computes the processing time of the of the $i$-th job of class $P_{k}$ which is executed. Constraint~\eqref{cst:startafter} (resp., \eqref{cst:startbefore}) ensures that a job which is scheduled after (before) another job (in accordance with the binary variable which defines the relative position of jobs) starts actually after (before) that job. Constraints~\eqref{cst:genduedatex1} and~\eqref{cst:genduedatex2} deal with the generalized due-date model again: in this case, with respect of class $P_{k}$, the constraints ensure that all jobs preceding the $i$-th are scheduled before job $i$ and all jobs following the $i$-th are scheduled after job $i$; the redundancy of constraints~\eqref{cst:genduedatedelta1}-\eqref{cst:genduedatedelta2} and~\eqref{cst:genduedatex1}-\eqref{cst:genduedatex2} in dealing with the generalized due-date model reduces the computational times to find an optimal solution. Constraints~\eqref{cst:xnocirculartwo} and~\eqref{cst:xnocircularthree} prevent the system to provide solutions in which some jobs constitute a closed/circular subsequence of executions. Constraint~\eqref{cst:relationshipdeltax} fixes the values of binary variables $x_{h,j,k,i}$ in accordance with the values of binary variables $\delta_{h,j,k,i}$; as a matter of fact, if the execution of the $i$-th job of class $P_{k}$ follows the execution of the $j$-th job of class $P_{h}$ ($\delta_{h,j,k,i} = 1$) then the $j$-th job of class $P_{h}$ definitely precedes the $i$-th job of class $P_{k}$ in the relative position of jobs ($x_{h,j,k,i} = 1$); such a constraint is equivalent to the simpler constraint $x_{h,j,k,i} \geq \delta_{h,j,k,i}$, but it provides better performance (lower computational times to solve the problem). Finally, constraints~\eqref{cst:defx}$\div$\eqref{cst:defLambda} define the type of the decision variables.

The proposed problem is a mixed-integer linear mathematical programming problem (MILP). The number of binary variables is $2(\sum_{k=1}^{K} N_{k})^{2}$.

\section{Test of the model} \label{sec:test}

Consider the numerical example for the single machine scheduling problem proposed in \cite{GiglioJOSH}, in which 4 jobs of class $P_{1}$ and 3 jobs of class $P_{2}$ must be executed. The data which characterize the problem are the following.
\begin{itemize}

\newpage
\item Due-dates and unitary tardiness costs

\vspace{6pt}
\begin{center}
\begin{tabular}{cc|cc}
\hline\noalign{\smallskip}
\multicolumn{2}{c|}{class $P_{1}$} & \multicolumn{2}{c}{class $P_{2}$}\\
\noalign{\smallskip}\hline\noalign{\smallskip}
$\alpha_{1,1} = 0.75$ & $dd_{1,1} = 19$ & $\alpha_{2,1} = 2$ & $dd_{2,1} = 21$\\
$\alpha_{1,2} = 0.5$ & $dd_{1,2} = 24$ & $\alpha_{2,2} = 1$ & $dd_{2,2} = 24$\\
$\alpha_{1,3} = 1.5$ & $dd_{1,3} = 29$ & $\alpha_{2,3} = 1$ & $dd_{2,3} = 38$\\
$\alpha_{1,4} = 0.5$ & $dd_{1,4} = 41$ & --- & ---\\
\noalign{\smallskip}\hline
\end{tabular}
\end{center}
\vspace{6pt}

\vspace{12pt}
\item Processing time bounds, unitary costs related to the deviation from the nominal processing times, and compression rates

\vspace{6pt}
\begin{center}
\begin{tabular}{ccc|ccc}
\hline\noalign{\smallskip}
\multicolumn{3}{c|}{class $P_{1}$} & \multicolumn{3}{c}{class $P_{2}$}\\
\noalign{\smallskip}\hline\noalign{\smallskip}
\multirow{2}{*}{$\beta_{1} = 1$} & $pt^{\mathrm{low}}_{1} = 4$ & \multirow{2}{*}{$\gamma_{1} = 1$} & \multirow{2}{*}{$\beta_{2} = 1.5$} & $pt^{\mathrm{low}}_{2} = 4$ & \multirow{2}{*}{$\gamma_{2} = 1$}\\
& $pt^{\mathrm{nom}}_{1} = 8$ & & & $pt^{\mathrm{nom}}_{2} = 6$ &\\
\noalign{\smallskip}\hline
\end{tabular}
\end{center}
\vspace{6pt}

\vspace{12pt}
\item Setup times and costs

\vspace{6pt}
\begin{center}
\begin{tabular}{cc}
\hline\noalign{\smallskip}
\multicolumn{2}{c}{times}\\
\noalign{\smallskip}\hline\noalign{\smallskip}
$st_{1,1} = 0$ & $st_{1,2} = 1$\\
$st_{2,1} = 0.5$ & $st_{2,2} = 0$\\
\noalign{\smallskip}\hline
\end{tabular}
\hspace{1cm}
\begin{tabular}{cc}
\hline\noalign{\smallskip}
\multicolumn{2}{c}{costs}\\
\noalign{\smallskip}\hline\noalign{\smallskip}
$sc_{1,1} = 0$ & $sc_{1,2} = 0.5$\\
$sc_{2,1} = 1$ & $sc_{2,2} = 0$\\
\noalign{\smallskip}\hline
\end{tabular}
\end{center}

\end{itemize}

\vspace{18pt}
By solving such an instance of the MILP problem, the following optimal solution is found.

\vspace{12pt}
\begin{center}
   \begin{tabular}{@{} cc|cccc|ccc @{}} 
      \toprule
      \multicolumn{2}{c|}{\multirow{2}{*}{$x_{h,j,k,i}$}} & \multicolumn{7}{c}{$(k,i)$}\\
      & & $1,1$ & $1,2$ & $1,3$ & $1,4$ & $2,1$ & $2,2$ & $2,3$\\
      \midrule 
	\multirow{7}{*}{$(h,j)$} & $1,1$ & 0 & 1 & 1 & 1 & 0 & 0 & 1\\
	& $1,2$ & 0 & 0 & 1 & 1 & 0 & 0 & 1\\
	& $1,3$ & 0 & 0 & 0 & 1 & 0 & 0 & 1\\
	& $1,4$ & 0 & 0 & 0 & 0 & 0 & 0 & 0\\
      \cmidrule{2-9} 
	& $2,1$ & 1 & 1 & 1 & 1 & 0 & 1 & 1\\
	& $2,2$ & 1 & 1 & 1 & 1 & 0 & 0 & 1\\
	& $2,3$ & 0 & 0 & 0 & 1 & 0 & 0 & 0\\
      \bottomrule
   \end{tabular}
\end{center}

\vspace{12pt}
\begin{center}
   \begin{tabular}{@{} cc|cccc|ccc @{}} 
      \toprule
      \multicolumn{2}{c|}{\multirow{2}{*}{$\delta_{h,j,k,i}$}} & \multicolumn{7}{c}{$(k,i)$}\\
      & & $1,1$ & $1,2$ & $1,3$ & $1,4$ & $2,1$ & $2,2$ & $2,3$\\
      \midrule 
	\multirow{7}{*}{$(h,j)$} & $1,1$ & 0 & 1 & 0 & 0 & 0 & 0 & 0\\
	& $1,2$ & 0 & 0 & 1 & 0 & 0 & 0 & 0\\
	& $1,3$ & 0 & 0 & 0 & 0 & 0 & 0 & 1\\
	& $1,4$ & 0 & 0 & 0 & 0 & 0 & 0 & 0\\
      \cmidrule{2-9} 
	& $2,1$ & 0 & 0 & 0 & 0 & 0 & 1 & 0\\
	& $2,2$ & 1 & 0 & 0 & 0 & 0 & 0 & 0\\
	& $2,3$ & 0 & 0 & 0 & 1 & 0 & 0 & 0\\
      \bottomrule
   \end{tabular}
\hspace{1cm}
   \begin{tabular}{@{} c|c @{}} 
      \toprule
      job $(k,i)$ & $u_{k,i}$\\
      \midrule 
	$1,1$ & 4\\
	$1,2$ & 0.5\\
	$1,3$ & 3\\
	$1,4$ & 0\\
      \midrule 
	$2,1$ & 0\\
	$2,2$ & 0\\
	$2,3$ & 0\\
      \bottomrule
   \end{tabular}
\end{center}

\newpage
With such values of the main decision variables, the optimal values of the other decision variables are

\vspace{6pt}
\begin{center}
   \begin{tabular}{@{} c|cccccc @{}} 
      \toprule
      job $(k,i)$ & $S_{k,i}$ & $pt_{k,i}$ & $T_{k,i}$ & $\Omega_{k,i}$ & $\Lambda_{k,i}$\\
      \midrule 
	$1,1$ & 12 & 4 & 0 & 1 & 0.5\\
	$1,2$ & 16.5 & 7.5 & 0 & 0 & 0\\
	$1,3$ & 24 & 5 & 0 & 0 & 0\\
	$1,4$ & 36 & 8 & 3.5 & 1 & 0.5\\
      \midrule 
	$2,1$ & 0 & 6 & 0 & 0 & 0\\
	$2,2$ & 6 & 6 & 0 & 0 & 0\\
	$2,3$ & 29 & 6 & 0 & 0.5 & 1\\
      \bottomrule
   \end{tabular}
\end{center}

\vspace{12pt}
The optimal job sequence provided by the solution is
\begin{equation*}
J_{2,1} \rightarrow J_{2,2} \rightarrow J_{1,1} \rightarrow J_{1,2} \rightarrow J_{1,3} \rightarrow J_{2,3} \rightarrow J_{1,4}
\end{equation*}
which is the same obtained with the methodology proposed in the paper in the nominal case. Moreover, the total cost of the optimal solution is 11.75, which corresponds to the value of the optimal cost-to-go in the initial state at $t_{0} = 0$, namely $J^{\circ}_{0,0,0}(0)$, provided in~\cite{GiglioJOSH}. However, it is worth noting that the optimal values of the decision variables (and thus the optimal job sequence) represent an open-loop solution that, when applied, maintains its validity as long as no disruption affects the nominal system behavior. Instead, the solution provided in~\cite{GiglioJOSH} is a closed-loop solution that is able to give the optimal decisions for any actual machine behavior.

\section{Experimental analysis} \label{sec:experiment}

An experimental analysis has been carried out to test the performance of the MILP model in solving the considered family scheduling problem. Randomly-generated instances of different sizes have been solved. All the parameters but $\gamma_{k}$ (which has been assumed deterministic with $\gamma_{k} = 1$ $\forall k = 1, \ldots, K$) have been considered stochastic variables uniformly distributed in the interval $[a,b]$, with $a$ and $b$ reported in the following table.

\vspace{6pt}
\begin{center}
   \begin{tabular}{@{} c|cc|c @{}} 
      \toprule
      \multirow{2}{*}{parameter} & \multicolumn{2}{|c|}{Uniform distribution $[a,b]$} & \multirow{2}{*}{note}\\
      & $a$ & $b$\\
      \midrule 
	$pt^{\mathrm{nom}}_{k}$ & 6 & 10 & \\
	$pt^{\mathrm{low}}_{k}$ & 2 & 6 & \\
	$dd_{k,i}$ & $dd_{k,i-1} + 0.5$ & $dd_{k,i-1} + 12$ & $dd_{k,0} = 10$\\
	$st_{h,k}$ & 1 & 3 & $st_{k,k} = 0$\\
	$sc_{h,k}$ & 0.5 & 2.5 & $sc_{k,k} = 0$\\
	$\alpha_{k,i}$ & 0.5 & 2.5 & \\
	$\beta_{k}$ & 0.5 & 2.5\\
      \bottomrule
   \end{tabular}
\end{center}

\vspace{6pt}
Only small instances of the problem have been taken into account, which can be solved in seconds or minutes. The results of such an analysis follow.

\vspace{12pt}
\fbox{2 classes -- 5 jobs per class} \hspace{.2cm} (Number of problems solved: 100)

\vspace{6pt}
\begin{center}
   \begin{tabular}{@{} c|c @{}} 
      \toprule
      \multicolumn{2}{c}{Problem size}\\
      \midrule 
	Binary variables & 200 \\
	Other variables & 61 \\
	Constraints & 1227 \\
      \bottomrule
   \end{tabular}
   \hspace{1cm}
   \begin{tabular}{@{} c|c @{}} 
      \toprule
      \multicolumn{2}{c}{Computational time (seconds)}\\
      \midrule 
	min & 0.14 \\
	max & 0.86 \\
	mean & 0.29 \\
	stdev & 0.09 \\
      \bottomrule
   \end{tabular}
\end{center}

\newpage
\fbox{2 classes -- 10 jobs per class} \hspace{.2cm} (Number of problems solved: 100)

\vspace{6pt}
\begin{center}
   \begin{tabular}{@{} c|c @{}} 
      \toprule
      \multicolumn{2}{c}{Problem size}\\
      \midrule 
	Binary variables & 800 \\
	Other variables & 121 \\
	Constraints & 8757 \\
      \bottomrule
   \end{tabular}
   \hspace{1cm}
   \begin{tabular}{@{} c|c @{}} 
      \toprule
      \multicolumn{2}{c}{Computational time (seconds)}\\
      \midrule 
	min & 8.17 \\
	max & 51.39 \\
	mean & 18.91 \\
	stdev & 7.69 \\
      \bottomrule
   \end{tabular}
\end{center}

\vspace{12pt}
\fbox{2 classes -- 15 jobs per class} \hspace{.2cm} (Number of problems solved: 20)

\vspace{6pt}
\begin{center}
   \begin{tabular}{@{} c|c @{}} 
      \toprule
      \multicolumn{2}{c}{Problem size}\\
      \midrule 
	Binary variables & 1800 \\
	Other variables & 181 \\
	Constraints & 28587 \\
      \bottomrule
   \end{tabular}
   \hspace{1cm}
   \begin{tabular}{@{} c|c @{}} 
      \toprule
      \multicolumn{2}{c}{Computational time (seconds)}\\
      \midrule 
	min & 181.39 \\
	max & 5013.48 \\
	mean & 1466.39 \\
	stdev & 1204.22 \\
      \bottomrule
   \end{tabular}
\end{center}

\vspace{12pt}
\fbox{3 classes -- 5 jobs per class} \hspace{.2cm} (Number of problems solved: 60)

\vspace{6pt}
\begin{center}
   \begin{tabular}{@{} c|c @{}} 
      \toprule
      \multicolumn{2}{c}{Problem size}\\
      \midrule 
	Binary variables & 450 \\
	Other variables & 91 \\
	Constraints & 3790 \\
      \bottomrule
   \end{tabular}
   \hspace{1cm}
   \begin{tabular}{@{} c|c @{}} 
      \toprule
      \multicolumn{2}{c}{Computational time (seconds)}\\
      \midrule 
	min & 16.36 \\
	max & 232.77 \\
	mean & 42.63 \\
	stdev & 29.03 \\
      \bottomrule
   \end{tabular}
\end{center}

\vspace{12pt}
Bigger instances of the scheduling problem take several hours to be solved.

\section{Comparison with the dynamic programming-based approach} \label{sec:comparison}

The performance of the dynamic programming-based (DP) approach proposed in~\cite{GiglioJOSH} and the performance of the MILP approach have been compared on some specific instances of the scheduling problem.
\begin{itemize}
\item DP approach $\rightarrow$ Implemented and solved with Matlab R2014a
\item MILP approach $\rightarrow$ Implemented and solved with Cplex 12.6 (with standard settings)
\end{itemize}

\vspace{6pt}
The results of such a comparison follow.

\vspace{12pt}
\fbox{2 classes -- 5 jobs per class} \hspace{.2cm}

\vspace{6pt}
\begin{center}
   \begin{tabular}{@{} c|c @{}} 
      \toprule
      & Time to find an optimal\\
      & solution (seconds)\\
      \midrule 
	DP approach & 0.52 \\
	MILP approach & 0.25 \\
      \bottomrule
   \end{tabular}
   \hspace{.75cm}
   \begin{tabular}{@{} c|c @{}} 
      \toprule
      \multicolumn{2}{c}{DP Problem size}\\
      \midrule 
	State space nodes & 61\\
      \bottomrule
   \end{tabular}
   \hspace{.75cm}
   \begin{tabular}{@{} c|c @{}} 
      \toprule
      \multicolumn{2}{c}{MILP Problem size}\\
      \midrule 
	Binary variables & 200 \\
	Other variables & 61 \\
	Constraints & 1227 \\
	MILP nodes & 527 \\
      \bottomrule
   \end{tabular}
\end{center}

\newpage
\fbox{2 classes -- 10 jobs per class} \hspace{.2cm}

\vspace{6pt}
\begin{center}
   \begin{tabular}{@{} c|c @{}} 
      \toprule
      & Time to find an optimal\\
      & solution (seconds)\\
      \midrule 
	DP approach & 1.47 \\
	MILP approach & 11.06 \\
      \bottomrule
   \end{tabular}
   \hspace{.75cm}
   \begin{tabular}{@{} c|c @{}} 
      \toprule
      \multicolumn{2}{c}{DP Problem size}\\
      \midrule 
	State space nodes & 221\\
      \bottomrule
   \end{tabular}
   \hspace{.75cm}
   \begin{tabular}{@{} c|c @{}} 
      \toprule
      \multicolumn{2}{c}{MILP Problem size}\\
      \midrule 
	Binary variables & 800 \\
	Other variables & 121 \\
	Constraints & 8757 \\
	MILP nodes & 13580 \\
      \bottomrule
   \end{tabular}
\end{center}

\vspace{6pt}
\fbox{2 classes -- 15 jobs per class} \hspace{.2cm}

\vspace{6pt}
\begin{center}
   \begin{tabular}{@{} c|c @{}} 
      \toprule
      & Time to find an optimal\\
      & solution (seconds)\\
      \midrule 
	DP approach & 3.17 \\
	MILP approach & 2903.38 \\
      \bottomrule
   \end{tabular}
   \hspace{.75cm}
   \begin{tabular}{@{} c|c @{}} 
      \toprule
      \multicolumn{2}{c}{DP Problem size}\\
      \midrule 
	State space nodes & 481\\
      \bottomrule
   \end{tabular}
   \hspace{.75cm}
   \begin{tabular}{@{} c|c @{}} 
      \toprule
      \multicolumn{2}{c}{MILP Problem size}\\
      \midrule 
	Binary variables & 1800 \\
	Other variables & 181 \\
	Constraints & 28587 \\
	MILP nodes & 2542617 \\
      \bottomrule
   \end{tabular}
\end{center}

\vspace{6pt}
\fbox{2 classes -- 20 jobs per class} \hspace{.2cm}

\vspace{6pt}
\begin{center}
   \begin{tabular}{@{} c|c @{}} 
      \toprule
      & Time to find an optimal\\
      & solution (seconds)\\
      \midrule 
	DP approach & 5.87 \\
	MILP approach & n.a. ($>$ 50000) \\
      \bottomrule
   \end{tabular}
   \hspace{.75cm}
   \begin{tabular}{@{} c|c @{}} 
      \toprule
      \multicolumn{2}{c}{DP Problem size}\\
      \midrule 
	State space nodes & 841\\
      \bottomrule
   \end{tabular}
   \hspace{.75cm}
   \begin{tabular}{@{} c|c @{}} 
      \toprule
      \multicolumn{2}{c}{MILP Problem size}\\
      \midrule 
	Binary variables & 3200 \\
	Other variables & 241 \\
	Constraints & 66717 \\
	MILP nodes & n.a. \\
      \bottomrule
   \end{tabular}
\end{center}

\vspace{6pt}
\fbox{3 classes -- 5 jobs per class} \hspace{.2cm}

\vspace{6pt}
\begin{center}
   \begin{tabular}{@{} c|c @{}} 
      \toprule
      & Time to find an optimal\\
      & solution (seconds)\\
      \midrule 
	DP approach & 4.99 \\
	MILP approach & 182.84 \\
      \bottomrule
   \end{tabular}
   \hspace{.75cm}
   \begin{tabular}{@{} c|c @{}} 
      \toprule
      \multicolumn{2}{c}{DP Problem size}\\
      \midrule 
	State space nodes & 541\\
      \bottomrule
   \end{tabular}
   \hspace{.75cm}
   \begin{tabular}{@{} c|c @{}} 
      \toprule
      \multicolumn{2}{c}{MILP Problem size}\\
      \midrule 
	Binary variables & 450 \\
	Other variables & 91 \\
	Constraints & 3790 \\
	MILP nodes & 492369 \\
      \bottomrule
   \end{tabular}
\end{center}

\vspace{6pt}
\fbox{3 classes -- 10 jobs per class} \hspace{.2cm}

\vspace{6pt}
\begin{center}
   \begin{tabular}{@{} c|c @{}} 
      \toprule
      & Time to find an optimal\\
      & solution (seconds)\\
      \midrule 
	DP approach & 54.69 \\
	MILP approach & n.a. ($>$ 50000) \\
      \bottomrule
   \end{tabular}
   \hspace{.75cm}
   \begin{tabular}{@{} c|c @{}} 
      \toprule
      \multicolumn{2}{c}{DP Problem size}\\
      \midrule 
	State space nodes & 3631\\
      \bottomrule
   \end{tabular}
   \hspace{.75cm}
   \begin{tabular}{@{} c|c @{}} 
      \toprule
      \multicolumn{2}{c}{MILP Problem size}\\
      \midrule 
	Binary variables & 1800 \\
	Other variables & 181 \\
	Constraints & 28435 \\
	MILP nodes & n.a. \\
      \bottomrule
   \end{tabular}
\end{center}

\vspace{6pt}
\fbox{4 classes -- 5 jobs per class} \hspace{.2cm}

\vspace{6pt}
\begin{center}
   \begin{tabular}{@{} c|c @{}} 
      \toprule
      & Time to find an optimal\\
      & solution (seconds)\\
      \midrule 
	DP approach & 91.31 \\
	MILP approach & n.a. ($>$ 50000) \\
      \bottomrule
   \end{tabular}
   \hspace{.75cm}
   \begin{tabular}{@{} c|c @{}} 
      \toprule
      \multicolumn{2}{c}{DP Problem size}\\
      \midrule 
	State space nodes & 4321\\
      \bottomrule
   \end{tabular}
   \hspace{.75cm}
   \begin{tabular}{@{} c|c @{}} 
      \toprule
      \multicolumn{2}{c}{MILP Problem size}\\
      \midrule 
	Binary variables & 800 \\
	Other variables & 121 \\
	Constraints & 8653 \\
	MILP nodes & n.a. \\
      \bottomrule
   \end{tabular}
\end{center}

\vspace{12pt}
The DP approach outperforms the MILP approach. Only in the case of very small instances of the problem, such as 2 classes and $\leq$ 5 jobs, the performance of the two approaches is comparable (both approaches solve the problem in less than 1 second), but in all the other cases the DP approach has significantly better performance than the MILP approach.

\section{Other mathematical programming models} \label{sec:model23}

Two different mixed-integer linear programming models are proposed in this section. The former is similar to the one proposed in subsection~\ref{sec:MIPmodel}, but it uses only the set of binary variables $\delta_{h,j,k,i}$; the latter exploits the stage-based state space representation, as introduced in \cite{GiglioJOSH}. However, these models do not show better performance (in terms of computational times to solve the considered instances of the scheduling problem) than that obtained with the model in subsection~\ref{sec:MIPmodel} of this technical report. 

\subsection{MILP model n.2}

\subsubsection{Variables and parameters}

Let, in addition to the variables and parameters introduced in subsection~\ref{sec:varpar}:
\begin{itemize}
\item $C_{k,i}$ : completion time of the $i$-th job of class $P_{k}$ (other decision variable)
\end{itemize}

\subsubsection{The model}

An alternative MILP model is the following.
\begin{equation}
\min \quad \sum_{k=1}^{K} \sum_{i=1}^{N_{k}} \alpha_{k,i} \, T_{k,i} + \sum_{k=1}^{K} \sum_{i=1}^{N_{k}} \beta_{k} \, \gamma_{k} \, u_{k,i} + \sum_{k=1}^{K} \sum_{i=1}^{N_{k}} \Omega_{k,i}
\end{equation}
subject to
\begin{equation} \label{cst:completionM2}
C_{k,i} = S_{k,i} + \Lambda_{k,i} + pt_{k,i} \qquad \forall k = 1, \ldots, K \quad \forall i = 1, \ldots, N_{k}
\end{equation}
\begin{equation} \label{cst:tardinessM2}
T_{k,i} \geq C_{k,i} - dd_{k,i} \qquad \forall k = 1, \ldots, K \quad \forall i = 1, \ldots, N_{k}
\end{equation}
\begin{equation} \label{cst:setupcostM2}
\Omega_{k,i} = \sum_{h=1}^{K} \sum_{j=1}^{N_{h}} sc_{h,k} \; \delta_{h,j,k,i} \qquad \forall k = 1, \ldots, K \quad \forall i = 1, \ldots, N_{k}
\end{equation}
\begin{equation} \label{cst:setuptimeM2}
\Lambda_{k,i} = \sum_{h=1}^{K} \sum_{j=1}^{N_{h}} st_{h,k} \; \delta_{h,j,k,i} \qquad \forall k = 1, \ldots, K \quad \forall i = 1, \ldots, N_{k}
\end{equation}
\begin{equation} \label{cst:alljobsM2}
\sum_{h=1}^{K} \sum_{j=1}^{N_{h}} \sum_{k=1}^{K} \sum_{i=1}^{N_{k}} \delta_{h,j,k,i} = \bigg( \sum_{k=1}^{K} N_{k} \bigg) - 1
\end{equation}
\begin{equation} \label{cst:onlyonepredecessorM2}
\sum_{h=1}^{K} \sum_{j=1}^{N_{h}} \delta_{h,j,k,i} \leq 1 \qquad \forall k = 1, \ldots, K \quad \forall i = 1, \ldots, N_{k}
\end{equation}
\begin{equation} \label{cst:onlyonesuccessorM2}
\sum_{k=1}^{K} \sum_{i=1}^{N_{k}} \delta_{h,j,k,i} \leq 1 \qquad \forall h = 1, \ldots, K \quad \forall j = 1, \ldots, N_{h}
\end{equation}
\begin{equation} \label{cst:genduedatedelta1M2}
\sum_{j=1}^{i} \delta_{k,i,k,j} = 0 \qquad \forall k = 1, \ldots, K \quad \forall i = 1, \ldots, N_{k}
\end{equation}
\begin{equation} \label{cst:genduedatedelta2M2}
\sum_{j=i+2}^{N_{k}} \delta_{k,i,k,j} = 0 \qquad \forall k = 1, \ldots, K \quad \forall i = 1, \ldots, N_{k}-2
\end{equation}
\begin{equation} \label{cst:consumptionM2}
u_{k,i} \leq \dfrac{pt^{\mathrm{nom}}_{k} - pt^{\mathrm{low}}_{k}}{\gamma_{k}} \qquad \forall k = 1, \ldots, K \quad \forall i = 1, \ldots, N_{k}
\end{equation}
\begin{equation} \label{cst:processingtimeM2}
pt_{k,i} = pt^{\mathrm{nom}}_{k} - \gamma_{k} u_{k,i} \qquad \forall k = 1, \ldots, K \quad \forall i = 1, \ldots, N_{k}
\end{equation}
\begin{multline} \label{cst:startafterM2}
S_{k,i} \geq C_{h,j} - M \, ( 1 - \delta_{h,j,k,i}) \\ \forall h = 1, \ldots, K \quad \forall j = 1, \ldots, N_{h} \quad \forall k = 1, \ldots, K \quad \forall i = 1, \ldots, N_{k} \quad (h,j) \neq (k,i)
\end{multline}
\begin{equation} \label{cst:firstjobM2}
C_{k,i} \geq pt_{k,i} - M \sum_{h=1}^{K} \sum_{j=1}^{N_{h}} \delta_{h,j,k,i} \qquad \forall k = 1, \ldots, K \quad \forall i = 1, \ldots, N_{k}
\end{equation}
\begin{equation} \label{cst:defdeltaM2}
\delta_{h,j,k,i} \in \{ 0 , 1 \} \qquad \forall h = 1, \ldots, K \quad \forall j = 1, \ldots, N_{h} \quad \forall k = 1, \ldots, K \quad \forall i = 1, \ldots, N_{k}
\end{equation}
\begin{equation} \label{cst:defuM2}
u_{k,i} \geq 0 \qquad \forall k = 1, \ldots, K \quad \forall i = 1, \ldots, N_{k}
\end{equation}
\begin{equation} \label{cst:defSM2}
S_{k,i} \geq 0 \qquad \forall k = 1, \ldots, K \quad \forall i = 1, \ldots, N_{k}
\end{equation}
\begin{equation} \label{cst:defCM2}
C_{k,i} \geq 0 \qquad \forall k = 1, \ldots, K \quad \forall i = 1, \ldots, N_{k}
\end{equation}
\begin{equation} \label{cst:defptM2}
pt_{k,i} \geq 0 \qquad \forall k = 1, \ldots, K \quad \forall i = 1, \ldots, N_{k}
\end{equation}
\begin{equation} \label{cst:defTM2}
T_{k,i} \geq 0 \qquad \forall k = 1, \ldots, K \quad \forall i = 1, \ldots, N_{k}
\end{equation}
\begin{equation} \label{cst:defOmegaM2}
\Omega_{k,i} \geq 0 \qquad \forall k = 1, \ldots, K \quad \forall i = 1, \ldots, N_{k}
\end{equation}
\begin{equation} \label{cst:defLambdaM2}
\Lambda_{k,i} \geq 0 \qquad \forall k = 1, \ldots, K \quad \forall i = 1, \ldots, N_{k}
\end{equation}

Constraint~\eqref{cst:completionM2} computes the completion time of the $i$-th job of class $P_{k}$ which is executed. Constraints~\eqref{cst:tardinessM2}$\div$\eqref{cst:processingtimeM2} have the same meanings of constraints~\eqref{cst:tardiness}$\div$\eqref{cst:processingtime} described in subsection~\ref{sec:MIPmodel}. Constraint~\eqref{cst:startafterM2} ensures that the $i$-th job of class $P_{k}$ is scheduled after the $j$-th job of class $P_{h}$ when $\delta_{h,j,k,i} = 1$. Constraint~\eqref{cst:firstjobM2} determines the completion time of the first job which is executed (since the first job is the only job for which it is $\sum_{h=1}^{K} \sum_{j=1}^{N_{h}} \delta_{h,j,k,i} = 0$). Finally, constraints~\eqref{cst:defdeltaM2}$\div$\eqref{cst:defLambdaM2} define the type of the decision variables.

The proposed problem is a mixed-integer linear mathematical programming problem (MILP). The number of binary variables is $(\sum_{k=1}^{K} N_{k})^{2}$.

\subsection{MILP model n.3}

The scheduling model can be represented through a state space model in which the system state, when a new decision has to be taken, namely at time instant $t_{j}$, $j = 0, 1, \ldots, N - 1$, with $N = \sum_{k = 1}^{K} N_{k}$, is the $(N+2)$-dimensional vector $\underline{x}_{j} = [n_{1,j}, \ldots, n_{K,j}, h_{j}, t_{j}]^{\mathrm{T}}$, being $n_{k,j}$, $k = 1, \ldots, K$, the number of jobs of class $P_{k}$ already completed at time instant $t_{j}$, and $h_{j}$ the class of the last served job. As a matter of fact, decision instants are discrete in time: they correspond to the initial time instant and the instants at which a job is completed, but for the last one.

\begin{figure*}[ht]
\centering
\psfrag{S0}[bc][Bl][.75][0]{$[0 \; 0 \; 0]^{T}$}
\psfrag{S1}[bc][Bl][.75][0]{$[1 \; 0 \; 1]^{T}$}
\psfrag{S2}[tc][Bl][.75][0]{$[0 \; 1 \; 2]^{T}$}
\psfrag{S3}[bc][Bl][.75][0]{$[2 \; 0 \; 1]^{T}$}
\psfrag{S4}[bc][Bl][.75][0]{$[1 \; 1 \; 1]^{T}$}
\psfrag{S5}[tc][Bl][.75][0]{$[1 \; 1 \; 2]^{T}$}
\psfrag{S6}[tc][Bl][.75][0]{$[0 \; 2 \; 2]^{T}$}
\psfrag{S7}[bc][Bl][.75][0]{$[3 \; 0 \; 1]^{T}$}
\psfrag{S8}[bc][Bl][.75][0]{$[2 \; 1 \; 1]^{T}$}
\psfrag{S9}[tc][Bl][.75][0]{$[2 \; 1 \; 2]^{T}$}
\psfrag{S10}[bc][Bl][.75][0]{$[1 \; 2 \; 1]^{T}$}
\psfrag{S11}[tc][Bl][.75][0]{$[1 \; 2 \; 2]^{T}$}
\psfrag{S12}[tc][Bl][.75][0]{$[0 \; 3 \; 2]^{T}$}
\psfrag{S13}[bc][Bl][.75][0]{$[4 \; 0 \; 1]^{T}$}
\psfrag{S14}[bc][Bl][.75][0]{$[3 \; 1 \; 1]^{T}$}
\psfrag{S15}[tc][Bl][.75][0]{$[3 \; 1 \; 2]^{T}$}
\psfrag{S16}[bc][Bl][.75][0]{$[2 \; 2 \; 1]^{T}$}
\psfrag{S17}[tc][Bl][.75][0]{$[2 \; 2 \; 2]^{T}$}
\psfrag{S18}[bc][Bl][.75][0]{$[1 \; 3 \; 1]^{T}$}
\psfrag{S19}[tc][Bl][.75][0]{$[1 \; 3 \; 2]^{T}$}
\psfrag{S20}[bc][Bl][.75][0]{$[4 \; 1 \; 1]^{T}$}
\psfrag{S21}[tc][Bl][.75][0]{$[4 \; 1 \; 2]^{T}$}
\psfrag{S22}[bc][Bl][.75][0]{$[3 \; 2 \; 1]^{T}$}
\psfrag{S23}[tc][Bl][.75][0]{$[3 \; 2 \; 2]^{T}$}
\psfrag{S24}[bc][Bl][.75][0]{$[2 \; 3 \; 1]^{T}$}
\psfrag{S25}[tc][Bl][.75][0]{$[2 \; 3 \; 2]^{T}$}
\psfrag{S26}[bc][Bl][.75][0]{$[4 \; 2 \; 1]^{T}$}
\psfrag{S27}[tc][Bl][.75][0]{$[4 \; 2 \; 2]^{T}$}
\psfrag{S28}[bc][Bl][.75][0]{$[3 \; 3 \; 1]^{T}$}
\psfrag{S29}[tc][Bl][.75][0]{$[3 \; 3 \; 2]^{T}$}
\psfrag{S30}[bc][Bl][.75][0]{$[4 \; 3 \; 1]^{T}$}
\psfrag{S31}[tc][Bl][.75][0]{$[4 \; 3 \; 2]^{T}$}
\psfrag{ST0}[bc][Bl][.8][0]{\it stage $0$}
\psfrag{ST1}[bc][Bl][.8][0]{\it stage $1$}
\psfrag{ST2}[bc][Bl][.8][0]{\it stage $2$}
\psfrag{ST3}[bc][Bl][.8][0]{\it stage $3$}
\psfrag{ST4}[bc][Bl][.8][0]{\it stage $4$}
\psfrag{ST5}[bc][Bl][.8][0]{\it stage $5$}
\psfrag{ST6}[bc][Bl][.8][0]{\it stage $6$}
\psfrag{ST7}[bc][Bl][.8][0]{\it stage $7$}
\includegraphics[scale=.4]{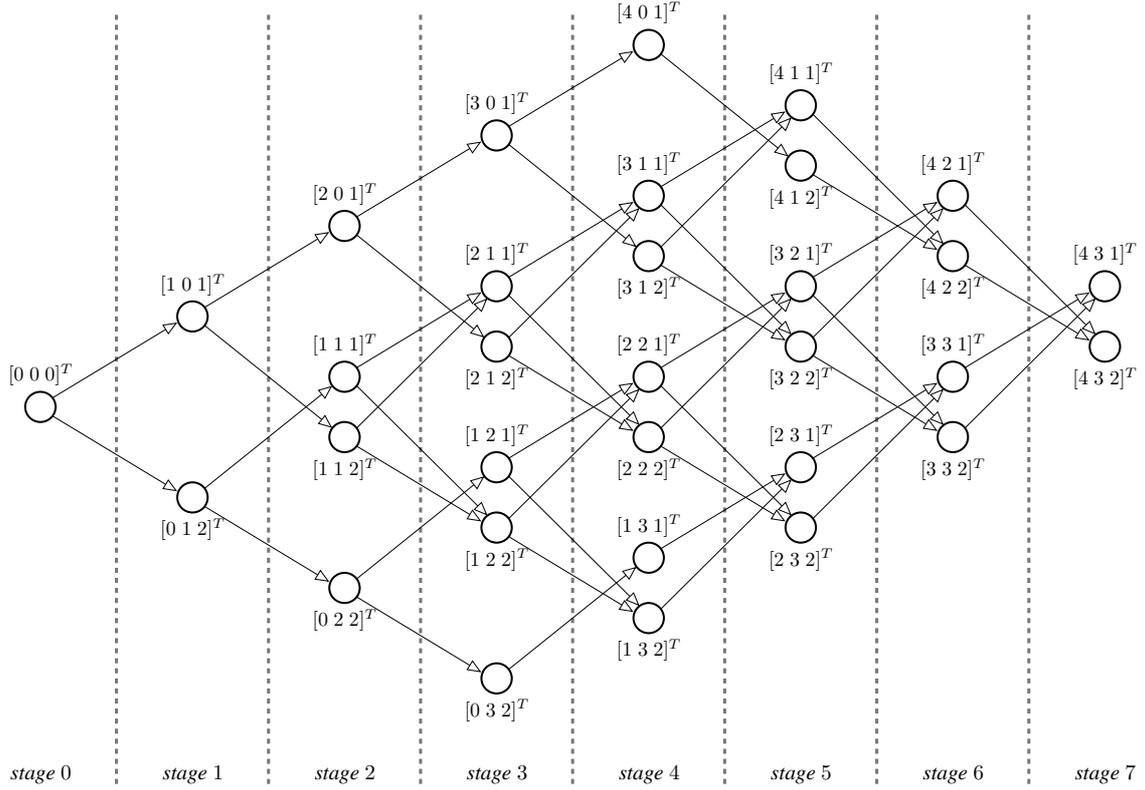}
\caption{State transition graph of the discrete part of the system state (in the case $K=2$, with $N_{1} = 4$, and $N_{2} = 3$)}
\label{fig:esS2_statediagram}
\end{figure*}

It is apparent that the system state does not change between two successive decision instants. The $(N+1)$-tuple $(n_{1,j}, \ldots, n_{K,j}, h_{j})$ will be referred to as the discrete part of the system state vector; it can be represented through a stage-based state transition graph in which the states at the $j$-th stage, $j \in \{ 0, 1, \ldots, N-1 \}$, are those such that $\sum_{k = 1}^{K} n_{k,j} = j$. An example of a state transition graph is reported in Fig.~\ref{fig:esS2_statediagram}, for the case of two classes of jobs, with $N_{1} = 4$ and $N_{2} = 3$.

The system state evolves with the application of the control $\underline{w}_{j} = \underline{f}(\underline{x}_{j})$ at the decision instants $t_{j}$, $j = 0, 1, \ldots, N - 1$. The control, which is a function of the system state, is relevant to the choice of the class of the next job to be served and the value of the processing time. Then, $\underline{w}_{j} = [ \delta_{1,j} , \ldots , \delta_{K,j} ,$ $\tau_{j} ]^{\mathrm{T}}$, where $\delta_{k,j} \in \{ 0, 1 \}$, $k = 1, \ldots, K$, is a binary decision variable whose value is $1$ if a job of class $P_{k}$ is selected for the $(j+1)$-th service, and $0$ otherwise, and $\tau_{j}$ is the processing time for the selected job. Obviously, it must be $\sum_{k=1}^{K} \delta_{k,j} = 1$ $\forall \, j = 0, 1, \ldots, N - 1$. Moreover, it it worth noting that the choice of considering the processing time $\tau_{j}$ as a decision variable instead of the amount of the continuous resource used to compress the processing time, namely $u_{j}$, does not modify the matter of the problem, as $u_{j} = \frac{1}{\sum_{k=1}^{K}{\gamma_{k} \, \delta_{k,j}}} \left( \sum_{k=1}^{K}{pt^{\mathrm{nom}}_{k} \, \delta_{k,j}} - \tau_{j} \right)$.

\vspace{6pt}

The state equations for the considered model are
\begin{equation} \label{equ:State_Equation}
\underline{x}_{j+1} = \left[ \begin{array}{c}
n_{1 , j+1} \\
\vdots\\
n_{K , j+1} \\
h_{j+1} \\
t_{j+1} \end{array} \right] = \left[ \begin{array}{c}
n_{1 , j} + \delta_{1 , j} \\
\vdots\\
n_{K , j} + \delta_{K , j} \\
\sum_{k=1}^{K}k \; \delta_{k,j} \\
t_{j} + \tau_{j} +\sum_{k=1}^{K} st_{h_{j},k} \; \delta_{k,j} \end{array} \right]
\end{equation}
and the initial state vector is $\underline{x}_{0} = [0, \ldots, 0, 0, 0]^{\mathrm{T}}$.

The objective is to find values $\delta_{k,j}$, $k = 1, \ldots, K$, and $\tau_{j}$, $j = 0, \ldots, N -1$, that minimize the objective function
\begin{equation} \label{equ:Objective_Function_j}
\sum_{j=0}^{N - 1} \Big[ \alpha_{h_{j+1},n_{h_{j+1},j+1}} \, \max{ \{ t_{j+1} - dd_{h_{j+1},n_{h_{j+1},j+1}}, 0 \}} + \beta_{h_{j+1}} \, ( pt^{\mathrm{nom}}_{h_{j+1}} - \tau_{j} ) + sc_{h_{j},h_{j+1}} \Big]
\end{equation}

\subsubsection{Variables and parameters}

Let, in addition to the variables and parameters introduced in section~\ref{sec:varpar}:
\begin{itemize}
\item $x_{k,i,j}$ : binary variable whose value is 1 if the $i$-th job of class $P_{k}$ is executed at the $j$-th stage, and $0$ otherwise (main decision variable). It is apparent that $\delta_{k,j} = \sum_{i=1}^{N_{k}} x_{k,i,j}$ (then variable $x_{k,i,j}$ can be employed in place of $\delta_{k,j}$).
\item $\tau_{j}$ : processing time of the job which is executed at the $j$-th stage (main decision variable)
\item $\tilde{\Omega}_{j}$ : setup cost which is paid at the $j$-th stage (other decision variable)
\item $\tilde{\Lambda}_{j}$ : setup time which is spent at the $j$-th stage (other decision variable)
\item $\tilde{S}_{j}$ : time instant at which the job which is executed at the $j$-th stage starts its execution (other decision variable)
\item $\tilde{C}_{j}$ : completion time of the job which is executed at the $j$-th stage (other decision variable)
\item $C_{k,i}$ : completion time of the $i$-th job of class $P_{k}$ (other decision variable)
\item $N$ : total number of jobs (parameter)
\end{itemize}

\subsubsection{The model}

An alternative MILP model is the following.
\begin{equation}
\min \quad \sum_{k=1}^{K} \sum_{i=1}^{N_{k}} \alpha_{k,i} \, T_{k,i} + \sum_{k=1}^{K} \sum_{i=1}^{N_{k}} \beta_{k} \big( pt^{\mathrm{nom}}_{k} - pt_{k,i} \big) + \sum_{j=0}^{N-1} \Omega_{j}
\end{equation}
subject to
\begin{equation} \label{cst:tardinessM3}
T_{k,i} \geq C_{k,i} - dd_{k,i} \qquad \forall k = 1, \ldots, K \quad \forall i = 1, \ldots, N_{k}
\end{equation}
\begin{equation} \label{cst:ptlowM3}
pt_{k,i} \geq pt^{\mathrm{low}}_{k} \qquad \forall k = 1, \ldots, K \quad \forall i = 1, \ldots, N_{k}
\end{equation}
\begin{equation} \label{cst:ptnomM3}
pt_{k,i} \leq pt^{\mathrm{nom}}_{k} \qquad \forall k = 1, \ldots, K \quad \forall i = 1, \ldots, N_{k}
\end{equation}
\begin{equation} \label{cst:setupcostM3}
\tilde{\Omega}_{j} \geq sc_{h,k} \, \bigg( \sum_{i=1}^{N_{h}} x_{h,i,j-1} + \sum_{i=1}^{N_{k}} x_{k,i,j} - 1 \bigg) \qquad \forall j = 1, \ldots, N-1 \quad \forall h,k = 1, \ldots, K
\end{equation}
\begin{equation} \label{cst:initialsetupcostM3}
\tilde{\Omega}_{0} = 0
\end{equation}
\begin{equation} \label{cst:setuptimeM3}
\tilde{\Lambda}_{j} \geq st_{h,k} \, \bigg( \sum_{i=1}^{N_{h}} x_{h,i,j-1} + \sum_{i=1}^{N_{k}} x_{k,i,j} -1 \bigg) \qquad \forall j = 1, \ldots, N-1 \quad \forall h,k = 1, \ldots, K
\end{equation}
\begin{equation} \label{cst:initialsetuptimeM3}
\tilde{\Lambda}_{0} = 0
\end{equation}
\begin{equation} \label{cst:sequenceM3}
\tilde{S}_{j} = \tilde{C}_{j-1} \qquad \forall j = 1, \ldots, N-1
\end{equation}
\begin{equation} \label{cst:initialM3}
\tilde{S}_{0} = 0
\end{equation}
\begin{equation} \label{cst:completionM3}
\tilde{C}_{j} = \tilde{S}_{j} + \tilde{\Lambda}_{j} + \tau_{j} \qquad \forall j = 0, \ldots, N-1
\end{equation}
\begin{equation} \label{cst:processingtimeM3}
\tau_{j} \geq pt_{k,i} - M \, ( 1 - x_{k,i,j}) \qquad \forall j = 0, \ldots, N-1 \quad \forall k = 1, \ldots, K \quad \forall i = 1, \ldots, N_{k}
\end{equation}
\begin{equation} \label{cst:startkiM3}
S_{k,i} \geq \tilde{S}_{j} - M \, ( 1 - x_{k,i,j}) \qquad \forall j = 0, \ldots, N-1 \quad \forall k = 1, \ldots, K \quad \forall i = 1, \ldots, N_{k}
\end{equation}
\begin{equation} \label{cst:completionkiM3}
C_{k,i} \geq \tilde{C}_{j} - M \, ( 1 - x_{k,i,j}) \qquad \forall j = 0, \ldots, N-1 \quad \forall k = 1, \ldots, K \quad \forall i = 1, \ldots, N_{k}
\end{equation}
\begin{equation} \label{cst:genduedateM3}
S_{k,i} \geq C_{k,i-1} \qquad \forall k = 1, \ldots, K \quad \forall i = 2, \ldots, N_{k}
\end{equation}
\begin{equation} \label{cst:onlyonceM3}
\sum_{k=1}^{K} \sum_{i=1}^{N_{k}} x_{k,i,j} = 1 \qquad \forall j = 0, \ldots, N-1
\end{equation}
\begin{equation} \label{cst:alljobsM3}
\sum_{j=0}^{N-1} \sum_{i=1}^{N_{k}} x_{k,i,j} = N_{k} \qquad \forall k = 1, \ldots, K
\end{equation}
\begin{equation} \label{cst:onlyoneM3}
\sum_{j=0}^{N-1} x_{k,i,j} = 1 \qquad \forall k = 1, \ldots, K \quad \forall i = 1, \ldots, N_{k}
\end{equation}
\begin{equation} \label{cst:defxM3}
x_{k,i,j} \in \{ 0 , 1 \} \qquad \forall k = 1, \ldots, K \quad \forall i = 1, \ldots, N_{k} \quad \forall j = 0, \ldots, N-1
\end{equation}
\begin{equation} \label{cst:deftauM3}
\tau_{j} \geq 0 \qquad \forall j = 0, \ldots, N-1
\end{equation}
\begin{equation} \label{cst:defSM3}
S_{k,i} \geq 0 \qquad \forall k = 1, \ldots, K \quad \forall i = 1, \ldots, N_{k}
\end{equation}
\begin{equation} \label{cst:defCM3}
C_{k,i} \geq 0 \qquad \forall k = 1, \ldots, K \quad \forall i = 1, \ldots, N_{k}
\end{equation}
\begin{equation} \label{cst:defptM3}
pt_{k,i} \geq 0 \qquad \forall k = 1, \ldots, K \quad \forall i = 1, \ldots, N_{k}
\end{equation}
\begin{equation} \label{cst:defTM3}
T_{k,i} \geq 0 \qquad \forall k = 1, \ldots, K \quad \forall i = 1, \ldots, N_{k}
\end{equation}
\begin{equation} \label{cst:defOmegajM3}
\tilde{\Omega}_{j} \geq 0 \qquad \forall j = 0, \ldots, N-1
\end{equation}
\begin{equation} \label{cst:defLambdajM3}
\tilde{\Lambda}_{j} \geq 0 \qquad \forall j = 0, \ldots, N-1
\end{equation}
\begin{equation} \label{cst:defSjM3}
\tilde{S}_{j} \geq 0 \qquad \forall j = 0, \ldots, N-1
\end{equation}
\begin{equation} \label{cst:defCjM3}
\tilde{C}_{j} \geq 0 \qquad \forall j = 0, \ldots, N-1
\end{equation}

Constraint~\eqref{cst:tardinessM3} computes, taking into account~\eqref{cst:defTM3} and the presence of $T_{k,i}$ in the cost function, the tardiness of the $i$-th job of class $P_{k}$ which is executed. Constraints~\eqref{cst:ptlowM3} and~\eqref{cst:ptnomM3} set the range of allowed values of the processing time of the $i$-th job of class $P_{k}$. Constraints~\eqref{cst:setupcostM3} and~\eqref{cst:setuptimeM3} compute, respectively, the setup cost and the setup time which are paid/spent at the $j$-th stage. Constraints~\eqref{cst:initialsetupcostM3} and~\eqref{cst:initialsetuptimeM3} set the initial setup cost and the initial setup time, respectively, to $0$. Since it can be shown that an optimal solution of the considered single machine family scheduling problem is characterized by the absence of idle times between the execution of two subsequent jobs (even if the cost function is not of regular type), constraint~\eqref{cst:sequenceM3} states that each job starts as soon as the previous one completes, and constraint~\eqref{cst:initialM3} sets the start time of the first job to $0$. Constraint~\eqref{cst:completionM3} computes the completion time of the job which is executed at the $j$-th stage. Constraint~\eqref{cst:processingtimeM3} sets the processing time of the job which is executed at the $j$-th stage, taking into account the values of binary variables $x_{k,i,j}$ as well as the range fixed by~\eqref{cst:ptlowM3} and~\eqref{cst:ptnomM3}. Constraints~\eqref{cst:startkiM3} and~\eqref{cst:completionkiM3} determine, respectively, the start time and the completion time of the $i$-th job of class $P_{k}$ which is executed. Constraint~\eqref{cst:genduedateM3} deals with the generalized due-date model as it ensures that, for what concerns class $P_{k}$, the $i$-th job is scheduled after the $(i-1)$-th job (that is, the $i$-th job starts after the completion time of the $(i-1)$-th job). Constraint~\eqref{cst:onlyonceM3} guarantees that a job is executed only once. Constraint~\eqref{cst:alljobsM3} ensures that all jobs are executed, whereas constraint~\eqref{cst:onlyoneM3} establishes that, in each stage, only one job can be served. Finally, constraints~\eqref{cst:defxM3}$\div$\eqref{cst:defCjM3} define the type of the decision variables.

The proposed problem is a mixed-integer linear mathematical programming problem (MILP). The number of binary variables is $(\sum_{k=1}^{K} N_{k})^{2}$.

\bibliographystyle{plain}
\bibliography{PaperScheduling}

\end{document}